\documentclass{amsart}
\usepackage{amsmath}
\usepackage{amssymb}
\usepackage{amsthm}
\usepackage{amscd}
\usepackage[all]{xy}

\newtheorem{Theo}{Theorem}
\newtheorem{Lem}[Theo]{Lemma}
\newtheorem{Prop}[Theo]{Proposition}
\newtheorem{Cor}[Theo]{Corollary}

\newtheorem{Prob}{Problem}

\newtheorem*{SM}{Super Multiplicity Theorem}

\newcommand{\Hom}{\mathrm{Hom}}
\newcommand{\de}{\mathrm{def}_p\,}
\newcommand{\deab}{\mathrm{def}_p^\mathrm{ab}\,}

\newcommand{\size}{\mathrm{size}}
\newcommand{\ab}{\mathrm{ab}}

\newcommand{\Z}{\mathbb{Z}}

\newcommand{\F}{\mathbb{F}}

\begin{document}
\title{On $p$-deficiency in groups}
\author[Y. Barnea]{Yiftach Barnea}
\author[J.-C. Schlage-Puchta]{Jan-Christoph Schlage-Puchta}
\begin{abstract}
Recently, Schlage-Puchta proved super multiplicity of $p$-deficiency for normal subgroups of $p$-power index. We extend this result to all normal subgroups of finite index. We then use the methods of the proof to show that some groups with non-positive $p$-deficiency have virtually positive $p$-deficiency. We also compute the $p$-deficiency in some cases such as Fuchsian groups and study related invariants: the lower and upper absolute $p$-homology gradients and the $p$-Euler characteristic.
\end{abstract}
\maketitle

\section{Introduction}
Let $\Gamma$ be a finitely generated group given by a presentation
$\Gamma$. We recall that the  \textit{deficiency} of $\langle X|R\rangle$ is $|X|-|R|$ which is denoted by $\mathrm{def} \langle X|R\rangle$ and the \textit{deficiency} of $\Gamma$, denoted by $\mathrm{def} \; \Gamma$, is the maximum of $\mathrm{def} \langle X|R\rangle$ over all possible finite presentations of $\Gamma$. If $\mathrm{def} \; \Gamma > 1$, then its abelianisation has more generators than relators, thus, it is infinite and in particular,  $\Gamma$ is infinite. However, groups $\Gamma$ with $\mathrm{def} \; \Gamma > 1$ are quite a small class of groups, for instance they cannot be torsion. Therefore, it is natural to look for a less restrictive criterion which still ensures that a group is infinite.

One such criterion is the Golod-Shafarevich inequality, see \cite{GS} and \cite{Golod}, in which relators are weighted according to their position in the Zassenhaus-filtration with respect to some prime number $p$. Using the Golod-Shafarevich inequality Golod \cite{Golod} was able to construct the first examples of finitely generated infinite $p$-groups and thus gave a negative answer to the General Burnside Problem. The most striking result concerning Golod-Shafarevich groups, i.e. groups in which the Golod-Shafarevich inequality holds, is Zelmanov's proof in \cite{Ze} that their pro-$p$ completion contains a non-abelian free pro-$p$ group.

Over the years, more examples of finitely generated infinite $p$-groups have been given. All of these constructions require quite difficult arguments. However, recently using remarkably simple arguments, more examples of finitely generated infinite $p$-groups were constructed in \cite{pdef}.

The first step was to generalize the notion of deficiency: Let $F$ be the free group over a finite set $X$, $w$ a non-trivial element of $F$. Define $\nu_p(w)$, the \textit{$p$-valuation} of $w$, to be the largest integer $k$, such that $w=v^{p^k}$ for some $v\in F$. Then define the \textit{$p$-deficiency} of $\langle X|R\rangle$, denoted by $\de\langle X|R\rangle$, to be $|X|-1-\sum_{r\in R}p^{-\nu_p(r)}$ and $\de \Gamma$, the \textit{$p$-deficiency} of $\Gamma$, to be the supremum of $\de\langle X|R\rangle$ taken over all presentations of $\Gamma$ with a finite generating set\footnote{We comment that it seems that the natural generalization of deficiency should be $|X|-\sum_{r\in R}p^{-\nu_p(r)}$ and indeed in \cite{BT} it is defined in this way. Nevertheless, because of the super multiplicity of our definition, see below, and connections to the $p$-homology we prefer to stick to it.}.

The key point of the construction was the proof of the super multiplicity of the $p$-deficiency for normal subgroups of $p$-power index, that is if $\Delta\triangleleft\Gamma$ is a normal subgroup of index $p$-power, then $\de\Delta\geq (\Gamma:\Delta)\; \de\Gamma$, see \cite[Theorem~2]{pdef}. In particular, if $\Gamma$ has positive (non-negative) $p$-deficiency, the same is true for all normal subgroups of $p$-power index in $\Gamma$. Notice that if $\Gamma$ has non-negative $p$-deficiency, then $\Gamma/([\Gamma,\Gamma]\Gamma^p)$ is non-trivial and therefore $\Gamma$ has a normal subgroup $\Delta$ of index $p$. Since $\Delta$ is finitely generated and has non-negative $p$-deficiency we have that $0<(\Delta:[\Delta,\Delta]\Delta^p)< \infty$. But $[\Delta,\Delta]\Delta^p$ is characteristic in $\Delta$ and therefore normal in $\Gamma$. Thus we can continue by induction and conclude that if $\Gamma$ has non-negative $p$-deficiency, then $\Gamma$ is infinite and furthermore has an infinite pro-$p$ completion.

Not only are groups with strictly positive $p$-deficiency infinite, they are big in other ways. Button and Thillaisundaram in \cite{BT}, and also Barnea, Ershov, Gonzales-Sanchez and Klopsch in unpublished work, showed that every group with positive $p$-deficiency is virtually Golod-Shafarevich (in most cases it is actually Golod-Shafarevich itself). Schlage-Puchta in \cite{pdef} showed that such a group has positive rank gradient. Moreover, it follows from the work of Lackenby in \cite{Lack1} that a finitely presented group $\Gamma$ of positive $p$-deficiency is large, that is, $\Gamma$ has a finite index subgroup that maps onto a non-abelian free group.

On the other hand, groups with zero $p$-deficiency, while infinite, do not have to be big. For instance, $\Z$ and $\Z^2$ have zero $p$-deficiency for all $p$ and $D_{\infty}$  and $$\langle x,y,z|x^2=y^4=z^4=xyz=1 \rangle$$ have zero $2$-deficiency and they are all not large.

Usually, knowing that the pro-$p$ completion of a group is non-trivial or even infinite tells us nothing about the pro-$p$ completion of a subgroup of finite index, unless the index is a $p$-power. Hence, our main theorem of this paper, which is a generalization of the super multiplicity to all normal subgroups of finite index, is somewhat surprising.
\begin{SM}
\label{SM}
Let $\Gamma$ be a finitely generated group and let $\Delta$ be a normal subgroup of finite index
in $\Gamma$. Then $\de\Delta\geq(\Gamma:\Delta)\de\Gamma$. In particular, if $\de\Gamma$ is non-negative (positive), then $\de\Delta$ is non-negative (positive), furthermore, $\Delta$ has a non-trivial pro-$p$ completion.
\end{SM}

Obviously this extends to subnormal subgroups of finite index. Nevertheless, we do not know the answer to the following problem.
\begin{Prob}
Let $\Gamma$ be a finitely generated group and let $\Delta$ be a subgroup of finite index
in $\Gamma$. If $\de\Gamma>0$, does it imply that $\de\Delta>0$?
\end{Prob}

Now, a group $\Gamma$ may contain a normal
subgroup $\Delta$ of finite index with $\de\Delta>0$ without $\Gamma$ having
positive $p$-deficiency. This is for example the case for all virtually free
groups which are generated by finite subgroups of $p'$-order because groups
of positive $p$-deficiency can be mapped onto a cyclic group of order $p$. The methods of the proof of the Super Multiplicity Theorem are flexible enough to enable us to study groups with zero $p$-deficiency. For instance, we prove that if $\Gamma\cong\langle X|R\rangle$, where $|X|$ is finite, $|R|>0$, $\de\langle X|R\rangle=0$, and each element in $R$ is a $q$-th power for some prime $q\neq p$, then $\Gamma$ contains a finite index subgroup of positive $p$-deficiency. Moreover, we can even study groups with negative $p$-deficiency. For example, we are able to show that the generalized triangle group $$\langle x,y|x^2=y^5=u(x,y)^5\rangle$$ which has negative $2$-deficiency has a normal subgroup of index $5$ of positive $2$-deficiency.

It seems impossible to give a general algorithm to compute the
$p$-deficiency of a group $\Gamma$, as doing so would involve all
presentations of $\Gamma$. However, if $\Gamma$ is restricted to some
well-understood class of groups things become easier. We compute
the $p$-deficiency for free products of cyclic groups and Fuchsian
groups of non-negative $p$-deficiency. It turns out that the obvious
presentations are the ones maximizing the $p$-deficiency, which provides evidence that the notion of $p$-deficiency is natural.

Finally we study a few new invariants of finitely generated groups:
Let $\Gamma$ be a finitely generated group. Write
$$d_p(\Gamma)=\dim_{\F_p}(\Gamma/([\Gamma,\Gamma]\Gamma^p))={ \dim_{\F_p} \Hom
(\Gamma, C_p)}.$$ Recall that Lackenby in \cite{Lack} defined the $p$-homology gradient of $\Gamma$ with respect to a particular chain of subgroups. Similarly we define the \textit{lower absolute $p$-homology gradient} of $\Gamma$
and the \textit{upper absolute $p$-homology gradient} of $\Gamma$
to be $$\alpha^{-}
(\Gamma)=\liminf \frac
{ d_p(\Delta)} {(\Gamma:\Delta)} \textrm{ and }  \alpha^{+}(\Gamma)=\limsup \frac
{ d_p(\Delta)} {(\Gamma:\Delta)},$$
respectively, where the limits are taken over all normal subgroups $\Delta$ of finite index in $\Gamma$. As an application of our methods we show that every finitely generated group of positive upper (lower) absolute
$p$-homology gradient maps onto a torsion group of positive upper (lower) absolute
$p$-homology gradient. We also define the \textit{$p$-Euler characteristic} of $\Gamma$:
\[
\chi_p(\Gamma) = -\sup_{(\Gamma:\Delta)<\infty}\frac{\de\Delta}{(\Gamma:\Delta)}.
\]
We then show that in this definition we only need to look at $\Delta$ which are normal in $\Gamma$. We also show that if $\Delta$ is a finite index subgroup of $\Gamma$, then $\chi_p(\Delta)=(\Gamma:\Delta) \chi_p(\Gamma)$. Further, we show that for virtually free groups and for Fuchsian groups the $p$-Euler characteristic coincides with the ordinary Euler characteristic and the hyperbolic volume respectively. We believe that answering the following problem will make an important contribution to geometric group theory.
\begin{Prob}
Find geometric interpretations of $p$-deficiency, upper and lower $p$-homology gradient and $p$-Euler characteristic.  \end{Prob}

\section{The $p$-size of a normal subgroup}

\setcounter{Theo}{0}

In this section we generalize the concept of $p$-deficiency
towards a relative notion measuring the size of a quotient of a
finitely generated group. Let $\Gamma$ be a group, $g\in \Gamma$ an
element of $\Gamma$. Then we define $\nu_{p, \Gamma}(g)$ to be the
supremum over all integers $n$ such that there exists some $h\in\Gamma$
with $g=h^{p^n}$. Let $\Gamma$ be a finitely generated group and let
$N\triangleleft\Gamma$ be a normal subgroup. Let $\{n_1, n_2, \ldots\} $ be
a set of elements generating $N$ as a normal subgroup. Then we call
$\sum_{j \geq 1} p^{-\nu_{p, \Gamma}(n_j)}$ the \textit{$p$-size} of
the generating system $\{n_1, n_2, \ldots\}$. We now define $\size_p(N,
\Gamma)$, the \textit{$p$-size of $N$ within $\Gamma$}, as the infimum
of the $p$-sizes of all normal generating systems of $N$. In this
section we study how $\size_p(N, \Gamma)$ behaves if we pass from
$\Gamma$ to a normal subgroup of finite index containing
$N$. The Super Multiplicity Theorem should convince the reader that this
notion both generalizes and improves the notion of $p$-deficiency.

Fix $p$ to be a prime number. For an integer $n$ we denote
$\nu_p(n)=\nu_{p, \mathbb{Z}}(n)$
to be the largest $k$ such that $p^k$ divides $n$. We prove the
following.
\begin{Lem}
\label{Lem:valuation}
Let $\Gamma$ be a group, $\Delta$ a normal subgroup of finite index $d$ in $\Gamma$, $g$ an element
of $\Delta$. Set $k=(C_\Gamma(g):C_{\Delta}(g))$. Then there exist elements $g_1, \ldots, g_{d/k}$, such that
\[
g^\Gamma=\bigcup_{i=1}^{d/k} g_i^{\Delta},
\]
and we have $\nu_{p, \Delta}(g)\geq\nu_{p, \Gamma}(g)-\nu_p(k)$.
\end{Lem}
\begin{proof}
Since $C_{\Gamma}(g)\Delta/\Delta \cong C_{\Gamma}(g)/\left(C_{\Gamma}(g) \cap \Delta\right)=C_\Gamma(g)/C_{\Delta}(g)$ we obtain that
$(C_{\Gamma}(g)\Delta:\Delta)=(C_\Gamma(g):C_{\Delta}(g))=k$ and thus, $$(\Gamma:C_{\Gamma}(g)\Delta)=\frac{(\Gamma:\Delta)}{(C_{\Gamma}(g)\Delta:\Delta)}=\frac{d}{k}.$$
Let $\{a_i\}_{i=1}^{d/k}$ be representatives of the right cosets of $C_{\Gamma}(g)\Delta$ in $\Gamma$ and let $\{c_j\}_{j=1}^{k}$ be representatives of the right cosets of $\Delta$ in $C_{\Gamma}(g)\Delta=\Delta C_{\Gamma}(g)$. With out loss of generality we can assume that $c_j \in C_{\Gamma}(g)$ for all $j$. As $\Delta$ is normal in $\Gamma$ we have that $\Delta c_ja_i=c_ja_i\Delta$ for all $i$ and $j$. Given $x \in \Gamma$, we have that $x=c_ja_in$ for some $i$ and $j$ and $n \in \Delta$. Since $c_j \in C_{\Gamma}(g)$ we have that $g^x=g^{c_ja_in}=g^{a_in}$. Write $g_i=g^{a_i}$ for all $i$. We conclude that
\[
g^\Gamma=\bigcup_{i=1}^{d/k} g_i^\Delta.
\]

For the inequality assume that $g=h^{p^m}$ for some $h\in \Gamma$, then $h\in C_
\Gamma(g)$. Notice that $C_{\Delta}(g)\triangleleft C_\Gamma(g)$. Hence, we obtain that both $h^{p^m} = g\in C_{\Delta}(g)$ and $h^
{k}$ are trivial in $C_\Gamma(g)/C_{\Delta}(g)$, and therefore $h^{(p^m, k)}$ is trivial
as well. As $(p^m, k)$ divides $\nu_p(k)$ this means that $h^{p^{\nu_p(k)}}$ is in $C_{\Delta}(g)<\Delta$. Thus, there exists an element
$h'=h^{p^{\nu_p(k)}}$ such that $g={h'}^{p^{m-\nu_p(k)}}$ and the inequality follows.
\end{proof}

\begin{Theo}
\label{thm:sizemult}
Let $\Gamma$ be a group, let $N\triangleleft\Gamma$ be a normal subgroup in $\Gamma$ and let
$\Delta \triangleleft\Gamma$ be a normal subgroup of finite index in $\Gamma$ containing $N$.
Then $$\mathrm{size}_{p} (N,\Delta)
\leq(\Gamma:\Delta)\mathrm{size}_{p}(N,\Gamma).$$
\end{Theo}
\begin{proof}
Suppose $\{n_1, n_2, \ldots\}$ is a set generating $N$ as a normal subgroup of $
\Gamma$ of $p$-size $\alpha$ with respect to $\Gamma$. Set $d=(\Gamma:\Delta)$ and $k_i=(C_\Gamma(n_i):C_\Delta
(n_i))$ for all $i$. As in Lemma~\ref{Lem:valuation} for each $i$ we write $$n_i^\Gamma = \bigcup_{j=1}^{d/
k_i} m_{ij}^\Delta.$$ Then $\{m_{ij} \}$ generates $N$ as a normal subgroup of $\Delta$.
The $p$-size of this set with respect to $\Delta$ is
\begin{multline}
\label{eq:multchain}
\sum_i \sum_{j=1}^{d/k_i} p^{-\nu_{p, \Delta}(m_{ij})} \leq \sum_i \frac{d}{k_i} p^{-\nu_{p, \Gamma}(n_i) +\nu_p(k_i)
} \leq  d\sum_i p^{-\nu_{p, \Gamma}(n_i)} = d\alpha
\end{multline}
Hence, there exists a generating set of $N$ as a normal subgroup of $\Delta$ with $p$-size at most $d\alpha=(\Gamma:\Delta)\alpha$. Our claim follows by taking the infimum over all generating sets of $N$ as a normal
subgroup of $\Gamma$.
\end{proof}

As a corollary we now
obtain the Super Multiplicity Theorem:
\begin{SM}
Let $\Gamma$ be a finitely generated group and let $\Delta$ be a normal subgroup of finite index
in $\Gamma$. Then $\de\Delta\geq(\Gamma:\Delta)\de\Gamma$. In particular, if $\de\Gamma$ is non-negative (positive), then $\de\Delta$ is non-negative (positive), furthermore, $\Delta$ has a non-trivial pro-$p$ completion.
\end{SM}
\begin{proof}
Consider a presentation $\Gamma=\langle X|R\rangle$.
Let $F$ be the free group with basis $X$ and let  $N$ be the normal subgroup in $F$ generated by
$R$. Taking $R$ as a generating set of $N$ we get that $$\size_{p}(N, F)
\leq |X|-1-\de\langle X|R\rangle.$$ Let $\overline{\Delta}$ be the pre-image of $\Delta$ in $F$. Since we have that
$(F:\overline{\Delta})=(\Gamma:\Delta)$ we obtain from Theorem~\ref{thm:sizemult} that $$\size_{p}(N,\overline
{\Delta})\leq (\Gamma:\Delta)\size_{p}(N,  F).$$ Let $Y$ be a basis of $\overline{\Delta}$ and recall that $|Y|-1=(F:\overline{\Delta})(|X|-1)=(\Gamma:\Delta)(|X|-1)$. We conclude that
\begin{eqnarray*}
\de\Delta & \geq & \underset{\langle Y|S\rangle\cong\Delta}{\sup_{S\subseteq \overline{\Delta}}}
|Y|-1-\sum_{s\in S} p^{-\nu_{p, \overline{\Delta}}(s)}\\
 & \geq & \underset{\langle S\rangle^{\overline{\Delta}} = N}{\sup_{S\subseteq \overline{\Delta}}}
|Y|-1-\sum_{s\in S} p^{-\nu_{p, \overline{\Delta}}(s)}\\
 & = & (\Gamma:\Delta)(|X|-1)-\size_{p}(N, \overline{\Delta})\\
 & \geq & (\Gamma:\Delta)(|X|-1)-(\Gamma:\Delta)\size_{p}(N, F)\\
  & \geq & (\Gamma:\Delta)\de\langle X|R\rangle.
\end{eqnarray*}
Our claim follows by taking the supremum over all presentations $\langle X|R\rangle$ .
\end{proof}

\section{Factoring out normal subgroups of small $p$-size}

Let $\Gamma$ be a finitely generated group and let $N$ be a normal subgroup in $\Gamma$. We
would like to show that if $\Gamma$ is large in some asymptotic way and if  $N$ has small $p$-size with respect to $\Gamma$, then $\Gamma/N$ is also large.

Let $\Gamma$ be a finitely generated group. Write
$$d_p(\Gamma)=\dim_{\F_p}(\Gamma/([\Gamma,\Gamma]\Gamma^p))={ \dim_{\F_p} \Hom
(\Gamma, C_p)}.$$ We recall from the introduction that the lower absolute $p$-homology gradient of $\Gamma$
and the upper absolute $p$-homology gradient of $\Gamma$
are $$\alpha^{-}
(\Gamma)=\liminf \frac
{ d_p(\Delta)} {(\Gamma:\Delta)} \textrm{ and }  \alpha^{+}(\Gamma)=\limsup \frac
{ d_p(\Delta)} {(\Gamma:\Delta)},$$
respectively, where the limits are taken over all normal subgroups $\Delta$ of finite index in $\Gamma$.

\begin{Theo}
\label{prop:psize}
Let $\Gamma$ be a finitely generated group and let $N$ be a normal subgroup in $\Gamma$. Then $\alpha^{-}
(\Gamma/N) \geq \alpha^{-}
(\Gamma) - \size_p(N,\Gamma)$ and $\alpha^{+}
(\Gamma/N) \geq \alpha^{+}
(\Gamma) - \size_p(N,\Gamma)$.
\end{Theo}
\begin{proof}
Let $\Delta$ be a normal subgroup of finite index in $\Gamma$ containing $N$. From Theorem~\ref
{thm:sizemult} we have that $\size_p(N, \Delta)\leq(\Gamma:\Delta)\size_p(N, \Gamma)$, thus we can
choose a set $\{n_1, n_2, \ldots\}$ which generates $N$ as a normal subgroup of $\Delta$ and satisfies  $$\sum p^{-\nu_{p, \Delta}(n_i)}\leq (\Gamma:\Delta)\size_p(N, \Gamma)+1.$$
Write $m=d_p(\Delta)$ and let $\varphi_1, \ldots, \varphi_{m}$ be a basis of the $\F_p$-vector space $\Hom(\Delta,
C_p)$. An element $\varphi\in\Hom(\Delta, C_p)$ defines a homomorphism $\overline{\varphi}:
\Delta/N\rightarrow C_p$ if and only if it maps each element of $N$ to 0, which is equivalent to the
statement that it maps each $n_i$ to 0. Since $\varphi$ is a homomorphism and $C_p$ has
exponent $p$ this condition is trivial whenever $n_i$ is a $p$-th power. Let $\{n_1, \ldots, n_\ell\}$
be a complete list of elements in the generating set which are not $p$-th powers of elements in $
\Delta$. Clearly $$\ell\leq \sum p^{-\nu_{p, \Delta}(n_i)}\leq (\Gamma:\Delta)\size_p(N, \Gamma)+1.$$
Express $\varphi$ as $\varphi=\sum_{i=1}^m x_i \varphi_i$. Then $\varphi$ defines a
homomorphism $\overline{\varphi}:\Delta/N\rightarrow C_p$ if and only if the coefficients $x_i$
satisfy the $\ell$ linear equations  $\sum_{i=1}^m x_i \varphi_i(n_j)=0$, $1\leq j\leq \ell$. The
dimension of the space of solutions of this system is at least the number of variables minus the
number of equations, thus,
\[
d_p(\Delta/N) \geq d_p(\Delta)-\ell \geq d_p(\Delta) - (\Gamma:\Delta)\size_p(N, \Gamma)-1
 \]
Now take a normal subgroup $\Delta'$ of $\Gamma/N$ of finite index. Then the pre-image of $
\Delta'$ under the canonical map $\Gamma\rightarrow\Gamma/N$ is a normal subgroup $\Delta$
of $\Gamma$ containing $N$ and satisfying $(\Gamma:\Delta)=(\Gamma/N:\Delta')$. Hence,
\[
\frac{d_p(\Delta')}{(\Gamma/N:\Delta')}  = \frac{d_p(\Delta/N)}{(\Gamma:\Delta)}\geq\frac{d_p
(\Delta)-1}{(\Gamma:\Delta)}-\size_p(N, \Gamma).
\]
Taking the $\liminf$ and $\limsup$ of this inequality implies our claims.
\end{proof}
As a corollary we obtain:
\begin{Cor}
Let $\Gamma$ be a finitely generated group of lower (upper) absolute $p$-homology gradient $\alpha>0$. Then for every
$\epsilon>0$ we have that $\Gamma$ maps surjectively onto a $p$-group with lower (upper) absolute $p$-homology
gradient $\geq\alpha-\epsilon$.
\end{Cor}
\begin{proof}
Let $k$ be an integer such that $\frac{p^{-k}}{1-p^{-1}}<\epsilon$ and let $g_1, g_2, g_3, \ldots$ be a list of all elements of $\Gamma$. Define $N$ as the normal subgroup generated by the elements $g_1^{p^k}, g_2^{p^{k+1}}, g_3^{p^{k+2}}, \ldots $. Then $\Gamma/N$ is a $p$-group, and $$\size_p(N, \Gamma)\leq \sum_{j=0}^\infty p^{-k-j} = \frac{p^{-k}}{1-p^{-1}}<\epsilon.$$ Hence, from the theorem we obtain that $\alpha^{-}(\Gamma/N)\geq\alpha^{-}(\Gamma)-\epsilon$ and $\alpha^{+}(\Gamma/N)\geq\alpha^{+}(\Gamma)-\epsilon$. Thus, $\Gamma/N$ is the desired image.
\end{proof}

\section{Groups of zero $p$-deficiency}

Let $\Gamma$ be a group given by a presentation $\langle X|R\rangle$ of zero
$p$-deficiency. Then $\Gamma$ needs not contain a subgroup of finite index of positive $p$-deficiency as seen in the examples in the introduction. In
this section we will consider under which conditions a group of zero $p$-deficiency has virtual positive $p$-deficiency.

\begin{Lem}
\label{Lem:power}
Let $\Gamma$ be a group given by a presentation $\langle X|R\rangle$, where
$|X|$ is finite, and $\de\langle X|R\rangle = 0$. Let $\Delta\triangleleft\Gamma$ be a normal
subgroup of index $n$ in $\Gamma$. Let $F$ be the free group on $X$ and let $\phi:F \to \Gamma$
be the map induced by the presentation $\langle X|R\rangle$. If there
exists an element $r\in R$, such that $r=g^e$ with $p \nmid e$ and $\phi(g) \not\in \Delta$, then $
\Delta$ has positive $p$-deficiency.
\end{Lem}
\begin{proof}
Let $N\triangleleft F$ be the normal subgroup generated by $R$ in $F$. By the Super Multiplicity Theorem $\Delta$ has non-negative $p$-deficiency. Inspecting the proof of the Super Multiplicity Theorem we see that $\Delta$ can only have zero $p$-deficiency, if $\size_p(N,\Delta)=(\Gamma:\Delta)\size_p(N, F)$. Going back to the proof of Theorem~\ref{thm:sizemult} we see that this is only possible if we have equality in (\ref{eq:multchain}). To have equality at the second position, one has to have $k_i=p^{\nu_p(k_i)}$ for all $i$, that is, all $k_i$ are powers of $p
$. In particular $(C_\Gamma(\phi(r)):C_\Delta(\phi(r)))$ is a power of $p$. But since $r=g^e$ we
have that $\phi(g)\in C_\Gamma(\phi(r))$, while $\phi(g)^e=\phi(r)=1$. In other words $C_\Gamma(\phi(r))/C_\Delta(\phi(r))$ is a $p
$-group containing an element of order dividing $e$, which is only possible if the order of this
element is 1, that is, $\phi(g)\in C_\Delta(\phi(r))\subseteq\Delta$. This contradicts the assumption $
\phi(g) \not\in \Delta$ and our claim is proven.
\end{proof}

For a word $w$ define the primitive $p'$-root of $w$ to be the
shortest word $v$, such that $v^n=w$ holds true for an integer $n$ not
divisible by $p$.
\begin{Theo}
\label{thm:Residual}
Let $\Gamma$ be a group given by a presentation $\langle X|R\rangle$, where $|X|$ is finite.
Suppose $\de\langle X|R\rangle=0$. Let $\langle X|S\rangle$ be the presentation obtained by
replacing every
word in $R$ by its primitive $p'$-root. Then $\Gamma$ contains a finite index subgroup of positive
$p$-deficiency or the kernel of the canonical map
$\Gamma\rightarrow\langle X|S\rangle$ is contained in the residual of
$\Gamma$.
\end{Theo}
\begin{proof}
We note that the kernel of the canonical map is $\left< S \right>^F/ \left< R \right>^F$. Hence, it
is enough to show that $s \left< R \right>^F$ is in the residual of
$\Gamma$ for every $s \in S$.

Suppose that $r\in R$, and $r=s^e$, where $p \nmid e$. In view of Lemma~\ref{Lem:power} it
suffices to show that
$\Gamma$ contains a normal subgroup $\Delta$ of finite index, such
that $s\left< R \right>^F \not\in\Delta$. Suppose this was not the case. Then $s\left<R \right>^F$ is
contained in every normal subgroup of finite index, hence, in the
residual of $\Gamma$.
\end{proof}
In particular, if $\Gamma$ is residually finite, then either the canonical map from $\langle X|R
\rangle$ to $\langle X|S\rangle$ is actually an isomorphism or $\Gamma$ contains a finite index
subgroup of positive $p$-deficiency. We call two groups \textit{finitely isomorphic} and write $G
\cong^*H$, if $G$ modulo its residual is isomorphic to $H$ modulo
its residual. We can now prove our claim from the introduction:
\begin{Theo}
\label{thm:Power}
Let $\langle X|R\rangle$ be a presentation, where $|X|$ is finite, $|R|>0$, and with zero $p$-deficiency. Suppose that each element in $R$ is a $q$-th power for some prime
$q\neq p$. Then $\Gamma\cong\langle
X|R\rangle$ contains a finite index subgroup of positive
$p$-deficiency.
\end{Theo}
\begin{proof} Let $S$ be as in the previous theorem. To apply the previous theorem we have to
show that $\langle X|R\rangle\not\cong^* \langle X|S\rangle$. It is enough to show that there exists a finite group $G$ such that the number of homomorphisms from $\langle X|R\rangle$ to $G$ is bigger than the number of homomorphisms from $\langle X|S\rangle$ to $G$ because the number of homomorphism to a finite group is
a property of a group modulo its residual.

We first note that for every finite group $G$ any homomorphism from $\langle X|S\rangle$ to $G$ induces a homomorphism from $\langle X|R\rangle$ to $G$ as $\langle X|S\rangle$ is a quotient of $\langle X|R\rangle$.
Furthermore, two distinct homomorphisms from $\langle X|S\rangle$ to $G$ will stay distinct when they
are induced to homomorphisms from $\langle X|R\rangle$ to $G$. Thus, the number of homomorphisms from $\langle X|R\rangle$ to $G$ is at least the number of homomorphisms from $\langle X|S\rangle$ to $G$.

We will find a homomorphism from $\langle X|R\rangle$ to $G$ which does not factor via a homomorphisms from $\langle X|S\rangle$ to $G$. This will finish the proof. Since free groups are residually $q$-groups we can choose a finite $q$-group $G$ and a homomorphism
$\alpha:F\rightarrow G$ such that some word $s_0\in S$ does not become trivial under $\alpha$.
Assume that $G$ is minimal with this property, that is, for every homomorphism $\beta:G\rightarrow
H
$ which is not injective we have $\beta\circ\alpha(s)=1$ for all $s \in S$. Since $G$ contains a non-trivial center, we can choose a subgroup $U$ in the center of order $q$. This subgroup is normal, and by
assumption we have that $\alpha(s)$ becomes trivial under the projection $G\rightarrow G/U$.
Hence,
$\alpha(s)\in U$ and for all $s\in S$ we have that $\alpha(s)$ has
order 1 or $q$, in particular, $\alpha(s_0)$ has order $q$. On the other hand, all elements in $R$
are mapped to the
trivial element under $\alpha$. We conclude that there exists a finite group $G$ and a homomorphism
$\alpha:F \rightarrow G$, which lifts to a homomorphism $\overline{\alpha}$ from $\langle X|R\rangle$ to $G$, but not to a
homomorphism from $\langle X|S\rangle$ to $G$.
\end{proof}

\begin{Cor}
\label{Cor:Large}
Let $\langle X|R\rangle$ be a presentation, where $|X|$ is finite, $R\neq\emptyset$, and with zero $p$-deficiency. Let $n \geq 2$ be an integer and $S=\left\{ r^n \mid r \in R \right\}$. Then $\Gamma_n=
\langle X|S\rangle$ contains a subgroup of finite index which has positive $p$-deficiency.
\end{Cor}
\begin{proof}
Clearly $\de\langle X|S\rangle \geq \de \langle X|R\rangle=0$ and if $p \mid n$ we have inequality.
Otherwise, take $q$ a prime such that $q \mid n$ and clearly every relation in $S$ is a $q$-power.
\end{proof}

The proofs of the results in this section were based on the fact that inequality~(\ref{eq:multchain}) is in general not sharp. It is hard to give general conditions under which there is an actual loss of a prescribed size, however, for an explicitly given presentation it is easy to compute finite images and check the kernels individually. As an example we prove the following.

\begin{Prop}
\label{Prop:negative}
Let $p, q$ be two distinct primes and $m$ a positive integer satisfying $m<q(1-\frac{1}{p})$. Let $F$ be a free group with basis $\left\{ x, y \right\}$ and let $w, v_1, \ldots, v_m \in F$. Let $\phi:F\rightarrow\F_q^2$ be the homomorphism defined by $\phi(x)=(1,0)$ and $\phi(y)=(0,1)$. Assume that none of the vectors $\phi(v_i)$ is a multiple of $\phi(w)$. Then the group $\Gamma=\langle x, y|w^p=v_1^q=\dots=v_m^q=1\rangle$ contains a normal subgroup of index $q$ and positive $p$-deficiency.
\end{Prop}
\begin{proof}
We first construct a homomorphism $\overline{\phi}:\Gamma\rightarrow C_q$, such that $w\in\ker\overline{\phi}$ and $v_i\not\in\ker\overline{\phi}$ for all $1\leq i\leq m$. If $\phi(w) \ne 0$, then let $\ell$ be the subspace of $\F_q^2$ spanned by $\phi(w)$. Since none of the vectors $\phi(v_i)$ is a multiple of $\phi(w)$, they are not contained in $\ell$. Otherwise, we recall that $\F_q^2$ contains $q+1$ distinct one-dimensional subspaces. As $m<q$, we can find  a one-dimensional subspace of $\F_q^2$ not containing any of the vectors $\phi(v_i)$, we call it $\ell$.

We now choose a non-trivial homomorphism $\rho:\F_q^2\rightarrow C_q$ mapping $\ell$ to zero. Let $\overline{\phi}=\rho\circ\phi$. From the definition of $\rho$, $\overline{\phi}(w^p)=0$ and $\overline{\phi}(v_i^q)=0$ since $C_q$ has exponent $q$. Thus, $\overline{\phi}$ maps the relations of $\Gamma$ to $0$ and, therefore, can be lifted to the required homomorphism. By a slight abuse of notation we call the lifted homomorphism from $\Gamma$ to $C_q$ also $\overline{\phi}$.

Let $\Delta$ be the kernel of $\overline{\phi}$ and $\overline{\Delta}$ the pre-image of $\Delta$ under the map $F\rightarrow\Gamma$ induced by the presentation $\Gamma=\langle x, y|w^p=v_1^q=\dots=v_m^q=1\rangle$. We would like to compute $(C_F(v_i^q):C_{\overline{\Delta}}(v_i^q))$. Clearly $v_i\in C_F(v_i^q)$, and by the construction of $\overline{\phi}$ we have $v_i\not\in\overline{\Delta}$, hence $(C_F(v_i^q):C_{\overline{\Delta}}(v_i^q))>1$. On the other hand $(F:\overline{\Delta})=q$, that is, $(C_F(v_i^q):C_{\overline{\Delta}}(v_i^q))$ divides $q$, we conclude that $(C_F(v_i^q):C_{\overline{\Delta}}(v_i^q))=q$.

Let $N$ be the normal subgroup generated by $\{w^p, v_1^q, \ldots, v_m^q\}$ in $F$. Then $N\leq\overline{\Delta}$, and we can estimate $\size_p(N, \overline{\Delta})$ using (\ref{eq:multchain}). We obtain that
\begin{multline*}
\size_p(N, \overline{\Delta})\leq \frac{qp^{-\nu_{p, F}(w^p)+\nu_p(q)}}{(C_F(w^p):C_{\overline{\Delta}}(w^p))}  + \sum_{i=1}^m \frac{qp^{-\nu_{p, F}(v_i^q)+\nu_p(q)}}{(C_F(v_i^q):C_{\overline{\Delta}}(v_i^q))}\\
 \leq q\cdot p^{-\nu_{p, F}(w^p)} + \sum_{i=1}^m p^{-\nu_{p, F}(v_i^q) }\leq \frac{q}{p} + m < q.
\end{multline*}
As $\overline{\Delta}$ is a subgroup of index $q$ in $F$ it is free with $q+1$ generators. We conclude that $\de\Delta\geq q+1-1-\size_p(N, \overline{\Delta})\geq q-m-\frac{q}{p}>0$, and our claim follows.
\end{proof}

The condition on $\phi$ looks quite technical, however, it is not very restrictive. To demonstrate this we give the following example.
\begin{Cor}
Let $u$ be a word in $x$ and $y$ for which the sum of the exponents of $y$ is not divisible by 5. Then the generalized triangle group $$\langle x,y|x^2=y^5=u(x,y)^5\rangle$$ has a normal subgroup of index 5 of positive 2-deficiency.
\end{Cor}
\begin{proof}
We take $p=2$, $q=5$, $m=2$ $w=x$, $v_1=y$, $v_2=u$ in the proposition. Clearly $2 < 5/2$. Now, $\phi(w)=(1,0)$, $\phi(v_1)=(0,1)$ and $\phi(u)=(a,b)$, where $a$ is the sum of the exponents of $x$ in $u$, and $b$ is the sum of the exponents of $y$ in $u$. Thus, $\phi(v_1)=(0,1)$ is not a multiple of $\phi(w)=(1,0)$ and $\phi(v_2)=(a,b)$ is not a multiple of $\phi(w)=(1,0)$ since $ b \not\equiv 0\pmod{5}$. Therefore, our group satisfies the conditions of the proposition.
\end{proof}
Notice that the abelianisation of $\Gamma=\langle x,y|x^2=y^5=u(x,y)^5\rangle$ is a subgroup of $C_2\times C_5$, hence, by Lemma~\ref{Lem:abelcompare} and Lemma~\ref{Lem:abelcompute} below these groups have $2$-deficiency at most $ -\frac{1}{2}$.


\section{Upper bounds for the $p$-deficiency via abelianisation}

In general it is difficult to compute the $p$-deficiency of a group. The $p$-deficiency is bounded above by the rank gradient and the lower absolute $p$-homology gradient, however, these invariants are also hard to compute, and often significantly differ from the $p$-deficiency, as the latter is much more sensitive to passing to subgroups.

In this section we define an abelianised version of the $p$-deficiency which can easily be computed and gives an upper bound for the $p$-deficiency. Moreover, this upper bound is tight in several cases.

Note first that the definition of $p$-deficiency $\Gamma$ depends on the notion of a presentation, which in turn depends on which category of groups we view $\Gamma$ as a member of. For example, the $p$-deficiency of pro-$p$-groups was defined in \cite{pdef} in the same way as for discrete groups. In fact, we could consider any category of groups for which free groups exist. However, it appears that only the case of abelian groups yields information which can be used for discrete groups.

We say that an abelian group $G$ has an abelian presentation $\langle X|R\rangle_{\ab}$, if $X$ is a set, $R$ is a set of elements of $\Z^X$, and $G\cong\Z^X/\langle R\rangle$. We define the abelian $p$-deficiency $\deab$ of $\langle X|R\rangle$ as $$|X|-1-\sum_{r\in R} p^{-\nu_{p, \Z^X}(r)}$$ and the abelian $p$-deficiency of $G$ as the supremum of $\deab\langle X|R\rangle$ taken over all abelian presentations of $G$.

\begin{Lem}
\label{Lem:abelcompare}
Let $\Gamma$ be a finitely generated group, $\Gamma^\ab=\Gamma/[\Gamma,\Gamma]$ the abelianisation of $\Gamma$. Then $\de\Gamma\leq\deab\Gamma^{ab}$.
\end{Lem}
\begin{proof}
Let $\langle X|R\rangle$ be a presentation of $\Gamma$ and let $F$ be the free group with a basis $X$. We will consider elements of $F$ as elements of $\Z^X$ by means of the canonical homomorphism $F\rightarrow\Z^X$. Then $\langle X|R\rangle$ is an abelian presentation of $\Gamma^\ab$. For each $r\in R$ we have $\nu_{p, F}(r)\leq\nu_{p, \Z^X}(r)$, hence, $\de\langle X|R\rangle\leq \deab\langle X|R\rangle$. Taking the supremum over all presentations $\langle X|R\rangle$ yields our claim.
\end{proof}

This simple observation is useful because the abelian $p$-deficiency can be computed easily.

\begin{Lem}
\label{Lem:abelcompute}
An abelian group $G\cong\oplus_{i=1}^s C_{e_i}\oplus \Z^r$ has abelian $p$-deficiency $r-1+\sum_{i=1}^s 1-p^{-\nu_p(e_i)}$.
\end{Lem}
\begin{proof}
The obvious presentation $\langle x_1, \ldots, x_s, y_1, \ldots, y_r|x_1^{e_1}=\dots=x_s^{e_s}=1\rangle$ has abelian $p$-deficiency $$r+s-1-\sum_{i=1}^s p^{-\nu_p(e_i)} = r-1+\sum_{i=1}^s 1-p^{-\nu_p(e_i)},$$ thus, $$\deab G \geq r-1+\sum_{i=1}^s 1-p^{-\nu_p(e_i)}.$$ It is left to show the inequality in the other direction. 

Let $\langle X|R\rangle$ be an abelian presentation of $G$. Let $n_j$ be the number of $r\in R$
with $\nu_{p,\mathrm{ab}}(r)=j$, and let $d_j$ be the number of indices $i$ such that $\nu_p(e_i)=j$. We now compute the number of cyclic factors of $G$ which have order divisible by $p^j$ for some integer $j$. From the explicit decomposition of $G$ we see that this number is $\sum_{i\geq j} d_i + r$. On the other hand, we have that $G\cong \Z^X/M_R\Z^R$, where $M_R$ is the $|X|\times|R|$-matrix with columns given by the elements of $R$. Let $S$ be the Smith normal form of $M_R$. Then $\Z^X/M_R\Z^R\cong \Z^X/S\Z^R$. All entries of $S$ are 0, with the possible exception of the entries $a_{ii}$, which satisfy $a_{11} \mid a_{22} \mid \dots$, hence, the number of cyclic factors of $G$ of order divisible by $p^j$ equals $|X|-\#\{i: p^j\nmid a_{ii}\}$. To bound this quantity note that the
number of rows, which do not vanish modulo $p^j$ in $S$ cannot be
larger than the number of rows of $M_R$ which do not vanish modulo $p^j$.
A relator $r$ gives rise to a non-vanishing row in $M_R$ if and only if $\nu_{p, \mathrm{ab}}(r) < j$, hence we get that
\[
\sum_{i\geq j} d_i + r \geq |X|-\sum_{i<j} n_i.
\]
We can perform the same computations starting from the standard presentation $\langle x_1, \ldots,
x_{s+r}|x_1^{e_1}=\dots=x_s^{e_s}=1\rangle$. Let $n_i'$ be the number of indices $j$ with $\nu_p
(e_j)=i$, $M'$ be the matrix obtained as before. As $M'$ is already diagonal the computation of the
Smith normal form reduces to a permutation of the diagonal elements and shuffling around $p'$-factors of the entries, that is, $M'$ and the Smith
normal form of $M'$ have the same number of rows vanishing modulo $p^j$ for every $j$. Hence,
\[
\sum_{i\geq j} d_i + s = r+s-\sum_{i<j} n_i'
\]
and comparing these two estimates we obtain that for all $j
$ $$\sum_{i<j} n_i'\leq (r+s-|X|)+\sum_{i<j} n_i.$$

From this we obtain our claim by a direct computation:
\begin{eqnarray*}
r+s-1-\sum_{i=1}^s p^{-\nu_p(e_i)} & =& r+s-1-\sum_{i\geq 0} n'_i p^{-i}\\
 & = & r+s-1-\sum_{i\geq 0} n'_i \big(1-\frac{1}{p}\big)\sum_{j\geq i} p^{-j}\\
 & = & r+s-1-\sum_{j\geq 0}\big(1-\frac{1}{p}\big)p^{-j} \sum_{i\leq j} n'_i\\
 & \geq & r+s-1-\sum_{j\geq 0}\big(1-\frac{1}{p}\big)p^{-j} \Big((s+r-|X|)+\sum_{i\leq j} n_i\Big)\\
 & = & r+s-1-(r+s-|X|)\sum_{j\geq 0}\big(1-\frac{1}{p}\big)p^{-j}\\
 && \qquad- \sum_{j\geq 0}\big(1-\frac{1}{p}\big)p^{-j}\sum_{i\leq j} n_i\\
  & = & |X|-1-\sum_{i\geq 0} n_i \big(1-\frac{1}{p}\big)\sum_{j\geq i} p^{-j}\\
 & = & |X|-1-\sum_{i\geq 0} n_i p^{-i}\\
 & = &\de\langle X|R\rangle.
\end{eqnarray*}
\end{proof}

\begin{Cor}
\label{Cor:freeproduct}
Let $e_1, \ldots, e_s$ be positive integers. Then the following free product  $\Gamma = C_{e_1}\ast\dots\ast C_{e_s}\ast F_r$ has $p$-deficiency $r-1+\sum_{i=1}^s 1-p^{-\nu_p(e_i)}$.
\end{Cor}
\begin{proof}
We have $\Gamma^\ab = C_{e_1}\oplus\dots\oplus C_{e_s}\oplus \Z^r$. Thus,
\begin{multline*}
r-1+\sum_{i=1}^s 1-p^{-\nu_p(e_i)} = \de\langle x_1, \ldots, x_s, y_1, \ldots, y_r|x_1^{e_1}=\dots=x_s^{e_s}=1\rangle\\
 \leq \de\Gamma \leq \deab\Gamma^\ab = r-1+\sum_{i=1}^s 1-p^{-\nu_p(e_i)},
\end{multline*}
hence, our claim follows.
\end{proof}

\section{$p$-deficiency of Fuchsian groups}

We can use a similar argument for Fuchsian groups.
Let $\Gamma$ be a finitely generated orientable Fuchsian group. If $\Gamma$
contains parabolic elements, then $\Gamma$ is a free product of
cyclic groups. This was considered above, so we can assume that $\Gamma$ is given by a
presentation
\[
\big\langle x_1, \ldots, x_r, u_1, v_1, \ldots, u_s, v_s\,\big|\\
x_1^{e_1}=\cdots=x_r^{e_r} = x_1\cdots x_r [u_1, v_1]\cdots [u_s, v_s]
= 1\big\rangle.
\]
We will refer to this presentation as the standard presentation
$\langle X|R\rangle$ of
$\Gamma$. The $p$-deficiency of this presentation is
$2s-2+\sum_{i=1}^r 1-p^{-\nu_p(e_i)}$. The abelianisation of $\Gamma$
consists of a free abelian factor of rank $2s$ and a finite factor,
which is obtained from $\bigoplus_{i=1}^r C_{e_i}$ by factoring out
the subgroup generated by $(1,\dots, 1)$, by which we mean the sum of
the generators of each cyclic factor. The $p$-component of this finite
group equals the $p$-component of $\bigoplus_{i=1}^r C_{e_i}$ with the
largest summand deleted, hence, if we assume that
$\nu_p(e_1)\geq\dots\geq\nu_p(e_r)$, we can apply Lemma~\ref{Lem:abelcompare} and Lemma~\ref{Lem:abelcompute} and find that
\begin{multline}
\label{eq:pdefbounds}
\de \langle X|R\rangle = 2s-2+\sum_{i=1}^r 1-p^{-\nu_p(e_i)}  \leq \\
\de\Gamma \leq 2s-1+\sum_{i=2}^r 1-p^{-\nu_p(e_i)}\leq \de \langle
X|R\rangle +1.
\end{multline}
We believe that the lower bound is correct. This believe is partly
founded on the naturalness of the standard presentation and partly on
the following Theorem.
\begin{Theo}
\label{Prop:Fuchs}
Suppose that one of the following statements holds true.
\begin{enumerate}
\item$s\geq 1$;
\item $p\geq 3$, and $e_1, e_2, e_3$ are divisible by $p$;
\item $p=2$, and $e_1, e_2, e_3, e_4$ are even;
\item $p=2$, $e_1$ and $e_2$ are divisible by 4, and $e_3$ is even.
\end{enumerate}
Then $\de\Gamma=\de\langle X|R\rangle$. If none of these conditions
holds true, then $\de\Gamma<0$.
\end{Theo}
For the proof we need the following theorem due to
Singerman~\cite{Sing}.
\begin{Theo}
Let $\Gamma$ be a Fuchsian group with elliptic generators $x_1,
\ldots, x_r$ of orders $e_1, \ldots, e_r$, respectively. Let $\Delta$
be a subgroup of index $n$ with associated coset action
$\varphi:\Gamma\rightarrow S_n$. For each $i$ let $\{c_{i1}, \ldots,
c_{im_i}\}$ be the possibly empty list of cycles of $\varphi(x_i)$,
which do not have length $e_i$. Then there is a bijection between
$\bigcup_{i=1}^r\{c_{i1}, \ldots, c_{im_i}\}$ and the elliptic
generators of $\Delta$, and this bijection maps a cycle $c_{ij}$ to
a generator of order $e_i/|c_{ij}|$.
\end{Theo}
We can now prove the theorem.
\begin{proof}
We start by showing that if $\Gamma$ is a group that satisfies  one of
the conditions (a)--(d), then $\Gamma$ has a normal subgroup $\Delta$ of
finite index which also satisfies one of these conditions and for
which the standard presentation $\langle Y|S\rangle$ of $\Delta$ has $p$-deficiency $\de\langle Y|S\rangle
=(\Gamma:\Delta)\de\langle X|R\rangle$. Assume first that $s\geq 1$. Then there exists an epimorphism $\phi:\Gamma\rightarrow C_2^2$ mapping all elliptic elements to the identity. Let $\Delta$ be the kernel of $\phi$. By Singerman's theorem $\Delta$ has 4 elliptic generators of order $e_1$, 4 generators of order $e_2$, \ldots. Let $s'$ be the number of commutators in the long relation of $S$, and let $\mu$ be the hyperbolic volume. Then we have
\[
2s'-2 + 4\sum_{i=1}^r 1-\frac{1}{e_i} = \mu(\Delta) = 4\mu(\Gamma) = 4\big(2s-2+\sum_{i=1}^r 1-\frac{1}{e_i} \big),
\]
that is, $s'=4s-3\geq 1$, and therefore, $\Delta$ satisfies (a). Moreover,
\[
\de\langle Y|S\rangle = 2s'-2+4\sum_{i=1}^r 1-p^{-\nu_p(e_i)} = 4\big(2s-2+\sum_{i=1}^r 1-p^{-\nu_p(e_i)}\big) = 4\de\langle X|R\rangle.
\]
Hence, $\Delta$ serves as the subgroup we are looking for.

In the cases (b)--(d) we have that $p$ divides $e_1$ and $e_2$. Hence we can construct a homomorphism $\phi:\Gamma\rightarrow C_p$ which maps $x_1$ to 1, $x_2$ to $-1$, and all other generators to 0. Let $\Delta$ be the kernel of $\phi$. By Singerman's theorem $\Delta$ has one elliptic generator of order $e_1/p$, one elliptic generator of order $e_2/p$ and $p$ elliptic generators of orders $e_3, \ldots, e_r$. Computing the hyperbolic volume
\[
2s'-2 + (1-\frac{p}{e_1})+(1-\frac{p}{e_2}) + p\sum_{i=3}^r 1-\frac{1}{e_i} = \mu(\Delta) = p\mu(\Gamma) = p\big(2s-2+\sum_{i=1}^r 1-\frac{1}{e_i} \big)
\]
we obtain $s'=ps$, and therefore,
\begin{multline*}
\de\langle Y|S\rangle = 2ps-2+(1-p^{1-\nu_p(e_1)})+(1-p^{1-\nu_p(e_2)}) + p\sum_{i=3}^r 1-p^{-\nu_p(e_i)}\\
 = p\big(2s-2+\sum_{i=1}^r 1-p^{-\nu_p(e_i)}\big) = p\de\langle X|R\rangle.
\end{multline*}
Hence, it suffices to check that $\Delta$ satisfies one of the conditions (b)--(d). If $\Gamma$ satisfies (b), then $p\mid e_3$ and $p\geq 3$. By Singerman's theorem $\Delta$ has $p$ elliptic generators of order $e_3$, hence $\Delta$ satisfies (b). If $\Gamma$ satisfies (c), then $e_3$ and $e_4$ are even, and $\Delta$ has 2 elliptic generators of order $e_3$, and two of order $e_4$, hence, $\Delta$ has four elliptic generators of even order, and therefore satisfies (c). If $\Gamma$ satisfies (d), then $\Delta$ has one elliptic generator of order $e_1/2$, one of order $e_2/2$, and two of order $e_3$. Since $e_1/2, e_2/2$ and $e_3$ are all even we see that $\Delta$ satisfies (c).

Next we construct a sequence of finite index subgroups $\Delta_i$ as follows. Put $\Delta_0=\Gamma$ and let $\Delta_{i+1}$ be a finite index normal subgroup in $\Delta_i$ satisfying one of the conditions (a)--(d), such that the $p$-deficiency of the standard presentation $\langle X_{i+1}|R_{i+1}\rangle$ of $\Delta_{i+1}$ satisfies $\de\langle X_{i+1}|R_{i+1}\rangle = (\Delta_i:\Delta_{i+1}) \de \langle X_i|R_i\rangle$. By the previous remark we can in fact obtain an infinite chain in this way. Each $\Delta_i$ is subnormal in $\Gamma$, hence, we have $\de\Delta_i\geq(\Gamma:\Delta_i)\de\Gamma$, putting this bound into (\ref{eq:pdefbounds}) we obtain that
\[
(\Gamma:\Delta_i)\de\Gamma\leq \de\Delta_i\leq \de\langle X_i| R_i\rangle + 1 = (\Gamma:\Delta_i)\de\langle X|R\rangle + 1.
\]
From this we have that $$\de\Gamma\leq\de\langle X|R\rangle + \frac{1}{(\Gamma:\Delta_i)}.$$ But $(\Gamma:\Delta_i)$ is unbounded, hence, $\de\Gamma\leq\de\langle X|R\rangle$, and in view of the obvious lower bound equality follows.

To prove the last statement of the theorem let $\Gamma$ be a group for which none of
the four criteria holds true. Then $s=0$, and from (\ref{eq:pdefbounds}) we obtain that $\de\Gamma\leq -1+\sum_{i=2}^r 1-p^{-\nu_p(e_i)}$. This is certainly negative, unless $r \geq 3$ and $p \mid e_3$, which we shall henceforth assume. But then we are in case (b), unless $p=2$. If $e_4$ is even, we are in case (c), and if $4\mid e_2$ we are in case (d). Hence it only remains to consider the case that $e_2$ and $e_3$ are even, but not divisible by 4. Let $\phi:\Gamma\rightarrow C_2$ be the homomorphism mapping $x_2$ and $x_3$ to 1 and all other generators to 0. Let $\Delta$ be the kernel of $\phi$. By Singerman's theorem $\Delta$ has one elliptic generator of order $e_2/2$, one of order $e_3/2$ and 2 of each of the orders $e_1, e_4, e_5, \ldots$. Since $e_2$ and $e_3$ are not divisible by 4 and $e_4, \ldots$ are odd, we have that $\Delta$ has exactly two elliptic generators of even order. By computing the hyperbolic volume as before we obtain that the long relation of $\Delta$ contains no commutators. From our previous discussion we find that $\Delta$ has negative 2-deficiency. But then $\Gamma$ must also have negative 2-deficiency and our claim follows.
\end{proof}

As an example consider the triangle group $\Delta(6,12,12)=\langle x,y|x^6=y^{12}=(xy)^
{12}\rangle$ and the generalized triangle group $\Gamma=\langle x,y|x^6=y^{12}=[x,y]^{12}\rangle$. It follows from the
previous computations that $\Delta(6,12,12)$ has  zero 2-deficiency and zero 3-deficiency, and $
\Gamma$ obviously has 2-deficiency and 3-deficiency $\geq 0$. Since all relators are sixth powers,
it follows from Theorem~\ref{thm:Power} that both groups contain finite index subgroups which
have positive 2-deficiency and 3-deficiency. From Lackenby \cite{Lack} we have that both groups
are large. For the triangle group this statement is an obvious consequence of the theory of
Fuchsian groups, however, it does not appear to be obvious for $\Gamma$.

\section{A $p$-Euler characteristics}
One common notion of the size of a group is its Euler characteristic. This invariant was originally
only defined for groups of finite homological dimension, but later on generalized to larger classes
of groups. For an overview we refer the reader to \cite{Chis}.
In \cite{pdef} the super multiplicity of the $p$-deficiency on normal subgroups of index $p$ was
used to define an Euler characteristics on pro-$p$-groups. Recall our definition from the introduction of the $p$-Euler characteristic for arbitrary finitely generated groups the following way:
\[
\chi_p(\Gamma) = -\sup_{(\Gamma:\Delta)<\infty}\frac{\de\Delta}{(\Gamma:\Delta)}.
\]
Having the superior Super Multiplicity Theorem
we can now prove:
\begin{Theo}
Let $\Gamma$ be a finitely generated group.
\begin{enumerate}
\item We have that
\[
\chi_p(\Gamma) = -\underset{\Delta\triangleleft\Gamma}{\sup_{(\Gamma:\Delta)<\infty}}\frac{\de\Delta}{(\Gamma:\Delta)}.
\]
\item If $\Delta$ is a finite index subgroup of $\Gamma$, then $\chi_p(\Delta)=(\Gamma:\Delta)
\chi_p(\Gamma)$.
\item If $\Gamma$ is virtually free, then $\chi_p(\Gamma)$ equals the ordinary Euler characteristic
of $\Gamma$.
\item If $\Gamma$ is a Fuchsian group, then $-\chi_p(\Gamma)$ equals the hyperbolic volume of $
\Gamma$.
\end{enumerate}
\end{Theo}
\begin{proof}
Part (a): The supremum taken over all normal subgroups is at most the supremum taken over all subgroups, thus, we have to show that the supremum taken over normal subgroups is large enough. Let $\Delta$ be a finite index subgroup of $\Gamma$. Then $N=\bigcap_{g\in\Gamma} \Delta^g$ is a normal subgroup of finite index in $\Gamma$. Applying Corollary 3 we obtain $\de N\geq (\Delta:N)\de\Delta$, and therefore
\[
\frac{\de N}{(\Gamma:N)}\geq \frac{(\Delta:N)\de\Delta}{(\Gamma:N)} = \frac{\de\Delta}{(\Gamma:\Delta)}.
\]
The supremum taken over all subgroups is therefore bounded by the supremum over all normal subgroups, and our claim follows.

Part (b): Since every finite index subgroup of $\Delta$ is of finite index in $\Gamma$ and the index is
multiplicative, the inequality $$-\chi_p(\Delta) \leq -(\Gamma:\Delta)\chi_p(\Gamma),$$ i.e $$\chi_p(\Delta) \geq (\Gamma:\Delta)\chi_p(\Gamma)$$ is trivial.

For the reverse inequality we use the fact that by (1) we may restrict the supremum to normal subgroups. Let $\Delta$ and $N$ be subgroups of finite index in $\Gamma$, where $N$ is normal in $\Gamma$.  Then $N\cap\Delta$ is normal in $\Delta$ and $$(\Gamma:\Delta)(\Delta:N\cap\Delta) = (\Gamma:N)(N:N\cap\Delta).$$ From the Super Multiplicity Theorem we obtain that
\[
\frac{\de(N\cap\Delta)}{(\Delta:N\cap\Delta)}\geq \frac{(N:N\cap\Delta)\de N}{(\Delta:N\cap\Delta)}
= \frac{(\Gamma:\Delta) \de N}{(\Gamma:N)}.
\]
Hence, if we take the supremum over all $N$ we obtain that $$-\chi_p(\Delta) \geq  -(\Gamma:\Delta)\chi_p(\Gamma),$$ i.e.  $$\chi_p(\Delta) \leq  (\Gamma:\Delta)\chi_p(\Gamma).$$

Parts (c) and (d): The Euler characteristics for virtually free groups and the hyperbolic volume for Fuchsian groups
are multiplicative on finite index subgroups. Hence, it suffices to
prove these statements for free groups and surface groups.
If $\Gamma$ is free on $d$ generators, then $\de\Gamma=d-1$, which is also its Euler characteristic. If $\Gamma$ is the fundamental group of an orientable surface of genus $g>0$, then $\de\Gamma=2g-2$, which is its hyperbolic volume.
\end{proof}

\bibliography{pdef_21_01_12}
\bibliographystyle{ijmart}

\end{document}